\documentclass[10pt]{amsart}
\usepackage{amsmath}
\usepackage{amssymb}
\usepackage{amsfonts}
\usepackage{longtable}
\usepackage{eucal}

\newtheorem{thm}{Theorem}

\newcommand{\ld}{\ldots}
\newcommand{\se}{\subseteq}
\newcommand{\mcd}{\!\cdot\!}
\renewcommand{\le}{\leqslant}
\renewcommand{\ge}{\geqslant}
\renewcommand{\S}{\mathfrak{S}}
\renewcommand{\SS}{\mathcal{S}}
\renewcommand{\AA}{\mathcal{A}}
\newcommand{\LL}{\mathcal{L}}
\newcommand\ord{\mathop\mathrm{ord}\nolimits}

\begin{document}
\renewcommand{\refname}{References}
\thispagestyle{empty}

\title{Finite simple groups with narrow prime spectrum}
\author{{Andrei V. Zavarnitsine}}%
\address{Andrei V. Zavarnitsine
\newline\hphantom{iii} Sobolev Institute of Mathematics,
\newline\hphantom{iii} pr. Koptyuga, 4,
\newline\hphantom{iii} 630090, Novosibirsk, Russia}
\email{zav@math.nsc.ru}%

\thanks{\sc A.V. Zavarnitsine,
Finite simple groups with narrow prime spectrum}
\thanks{\rm Supported by RFBR, grants 06-01-39001 and 08-01-00322;
by the Council of the President (project NSh-344.2008.1), by the Russian Science Support
Foundation (grant of the year 2008); and by SB RAS, Integration Project 2006.1.2}

\maketitle
{\small
\begin{quote}
\noindent{\sc Abstract. }
We find the nonabelian finite simple groups with order prime divisors not exceeding $1000$.
More generally, we determine the sets of nonabelian finite simple groups whose maximal order prime
divisor is a fixed prime less than $1000$.
Our results are based on calculations in the computer algebra system GAP.
 \end{quote}
}

\section{Introduction}

In the study of arithmetical properties of finite simple groups, it is often desirable to know
the groups whose order prime divisors are not too big. Given a finite group $G$, we denote by
$\pi(G)$ the {\em prime spectrum\/} of $G$, i.~e. the set of prime divisors of the order
$|G|$. Let $\pi$ be a finite set of primes. Since the orders of the simple groups are known,
it is possible, in principle, to determine all groups $G$ with $\pi(G)=\pi$.  Although one can
find these groups by hand, computer help usually facilitates the task and is less error-prone.
The main result of this paper is the following assertion:

\begin{thm}\label{mt} There are a total of\/ $1972$ finite nonabelian simple groups $G$ all of whose
order prime divisors do not exceed\/ $1000$.
\end{thm}

The groups from theorem \ref{mt} are listed in Tables
\ref{tab}--\ref{alph} of this paper. In naming simple groups, we mostly use the notation from \cite{atl}.

\section{The algorithm}

It was observed in \cite[p.\,51]{m94} that, given a finite set of primes $\pi$, there are only
finitely many finite (nonabelian) simple groups $G$ such that $\pi(G)\se \pi$. The gist of the
algorithm for finding these groups was also outlined in the proof of
\cite[Lemma 2]{m94}. One should have $G\in\SS\cup\AA_\pi\cup\LL_\pi$, where $\SS$ is the set of $26$
sporadic groups, $\AA_\pi$ is the set of the alternating groups $A_n$ with $n\le p-1$, where
$p$ is the smallest prime greater that all the primes in $\pi$, and $\LL_\pi$ is the set of
groups of Lie type of rank $l$ defined over a field of order $q=p^k$, where $p\in \pi$, $k\le
t$, $l\le \max\{8,t\}$, and
$$t=\max_{r\in \pi\setminus\{p\}}\ord_r p.$$
Hence, finding all such groups $G$ reduces to checking the orders of finitely many groups. An
implementation of this algorithm (with some technical improvements) was performed in the
computer algebra system GAP \cite{gap}. Theorem \ref{mt} is the result of running the program
for $\pi=\{2,3,5,\ld,997\}$.

\section{The tables}

Given a prime $p$, we denote by $\S_p$ the set of nonabelian finite simple groups $G$ such
that $p\in\pi(G)\se\{2,3,5,\ld,p\}$. It is clear that the set of {\em all\/} nonabelian finite
simple groups is the disjoint union of the finite sets $\S_p$ for all primes $p$. Observe that
the first two sets $\S_2$ and $\S_3$ are trivially empty, whereas the sets $\S_p$, $p\ge 5$
are always nonempty because they contain some generic elements (see below). The groups from
Theorem~\ref{mt} thus constitute $\S_5\cup\ld\cup\S_{997}$.

In Table \ref{tab}, we list the elements of $\S_p$ for $p\le 100$.  The upper bound of $100$
was chosen so as to include all the sporadic groups. The members of $\S_p$ are ordered
according to the size of their prime spectrum (same for Table \ref{ng}). The number of groups
in each set $\S_p$ is given after the symbol "\#". For each group, we also give the prime
decomposition of the order (except for $A_n$, $n\ge 23$, whose order decompositions occupy too
much space and are not too interesting).

The sets $\S_p$ for larger primes $p$ are described in a more abbreviated form. Denote by $p'$
the smallest prime greater than $p$. Then $\S_p$, $p\ge 5$, always contains the groups
$$L_2(p),A_p,A_{p+1},\ld,A_{p'-1}$$ which we call the {\em generic\/} elements of $\S_p$. The
primes $p\in\{100,\ld,1000\}$ such that $\S_p$ contains no simple groups other than the
generic ones are listed in Table~\ref{gen}. In this case, we have $|\S_p|=p'-p+1$. Note that
the only such primes $p$ less than $100$ are $59$ and $89$ as follows from Table~\ref{tab}.
The non-generic elements of $\S_p$ for the primes $p$ not listed in Table~\ref{gen} are given
in Table~\ref{ng}.

Finally, in Table \ref{alph} we list all the groups from Theorem \ref{mt} collected in series
so that the membership of a given simple group could be easily checked. Isomorphic groups like
$L_2(7)\cong L_3(2)$ are all included. The parameters $n,k$ are positive integers, $p$ is a
prime, and $q$ is a prime power. Long sequences of consecutive primes or integers are
abbreviated with an ellipsis.


{\footnotesize
\begin{longtable}[h]{l|l}
\caption {\bf Nonabelian simple groups $G$ with  $\pi(G)\se\{1,2,\ld,100\}$}\label{tab}\\
\hline
\multicolumn{1}{c|}{$^{\vphantom{A^A}}G$} & \multicolumn{1}{c}{$|G|$}\\
\hline
\endfirsthead
\multicolumn{2}{c}{$\S_5$ \hfill $^{\vphantom{A^A}}\pi(G)\se\{2,3,5\}$ \hfill $\#\,3$}  \\
\hline
 $^{\vphantom{A^A}}$$A_{5}\cong L_2(4)\cong L_2(5)$ & $2^2 \mcd 3 \mcd 5$\\
$A_{6}\cong L_2(9)$ & $2^3 \mcd 3^2 \mcd 5$\\
$S_{4}(3)\cong U_4(2)$ & $2^6 \mcd 3^4 \mcd 5$\\
\hline
\multicolumn{2}{c}{$\S_{7}$ \hfill $^{\vphantom{A^A}}7\in\pi(G)\se\{2,3,5,7\}$ \hfill $\#\,15$}  \\
\hline
 $^{\vphantom{A^A}}$$L_{2}(7)\cong L_3(2)$ & $2^3 \mcd 3 \mcd 7$\\
$L_{2}(8)$ & $2^3 \mcd 3^2 \mcd 7$\\
$U_{3}(3)$ & $2^5 \mcd 3^3 \mcd 7$\\
$A_{7}$ & $2^3 \mcd 3^2 \mcd 5 \mcd 7$\\
$L_{2}(49)$ & $2^4 \mcd 3 \mcd 5^2 \mcd 7^2$\\
$U_{3}(5)$ & $2^4 \mcd 3^2 \mcd 5^3 \mcd 7$\\
$L_{3}(4)$ & $2^6 \mcd 3^2 \mcd 5 \mcd 7$\\
$A_{8}\cong L_4(2)$ & $2^6 \mcd 3^2 \mcd 5 \mcd 7$\\
$A_{9}$ & $2^6 \mcd 3^4 \mcd 5 \mcd 7$\\
$J_2$ & $2^7 \mcd 3^3 \mcd 5^2 \mcd 7$\\
$A_{10}$ & $2^7 \mcd 3^4 \mcd 5^2 \mcd 7$\\
$U_{4}(3)$ & $2^7 \mcd 3^6 \mcd 5 \mcd 7$\\
$S_{4}(7)$ & $2^8 \mcd 3^2 \mcd 5^2 \mcd 7^4$\\
$S_{6}(2)$ & $2^9 \mcd 3^4 \mcd 5 \mcd 7$\\
$O^+_{8}(2)$ & $2^{12} \mcd 3^5 \mcd 5^2 \mcd 7$\\
\hline
\multicolumn{2}{c}{$\S_{11}$ \hfill $^{\vphantom{A^A}}11\in\pi(G)\se\{2,3,5,7,11\}$ \hfill $\#\,10$}  \\
\hline
 $^{\vphantom{A^A}}$$L_{2}(11)$ & $2^2 \mcd 3 \mcd 5 \mcd 11$\\
$M_{11}$ & $2^4 \mcd 3^2 \mcd 5 \mcd 11$\\
$M_{12}$ & $2^6 \mcd 3^3 \mcd 5 \mcd 11$\\
$U_{5}(2)$ & $2^{10} \mcd 3^5 \mcd 5 \mcd 11$\\
$M_{22}$ & $2^7 \mcd 3^2 \mcd 5 \mcd 7 \mcd 11$\\
$A_{11}$ & $2^7 \mcd 3^4 \mcd 5^2 \mcd 7 \mcd 11$\\
$McL$ & $2^7 \mcd 3^6 \mcd 5^3 \mcd 7 \mcd 11$\\
$HS$ & $2^9 \mcd 3^2 \mcd 5^3 \mcd 7 \mcd 11$\\
$A_{12}$ & $2^9 \mcd 3^5 \mcd 5^2 \mcd 7 \mcd 11$\\
$U_{6}(2)$ & $2^{15} \mcd 3^6 \mcd 5 \mcd 7 \mcd 11$\\
\hline
\multicolumn{2}{c}{$\S_{13}$ \hfill $^{\vphantom{A^A}}13\in\pi(G)\se\{2,3,5,\ld,13\}$ \hfill $\#\,27$}  \\
\hline
 $^{\vphantom{A^A}}$$L_{3}(3)$ & $2^4 \mcd 3^3 \mcd 13$\\
$L_{2}(25)$ & $2^3 \mcd 3 \mcd 5^2 \mcd 13$\\
$U_{3}(4)$ & $2^6 \mcd 3 \mcd 5^2 \mcd 13$\\
$S_{4}(5)$ & $2^6 \mcd 3^2 \mcd 5^4 \mcd 13$\\
$L_{4}(3)$ & $2^7 \mcd 3^6 \mcd 5 \mcd 13$\\
$^2F_4(2)'$ & $2^{11} \mcd 3^3 \mcd 5^2 \mcd 13$\\
$L_{2}(13)$ & $2^2 \mcd 3 \mcd 7 \mcd 13$\\
$L_{2}(27)$ & $2^2 \mcd 3^3 \mcd 7 \mcd 13$\\
$G_2(3)$ & $2^6 \mcd 3^6 \mcd 7 \mcd 13$\\
$^3D_4(2)$ & $2^{12} \mcd 3^4 \mcd 7^2 \mcd 13$\\
$Sz(8)$ & $2^6 \mcd 5 \mcd 7 \mcd 13$\\
$L_{2}(64)$ & $2^6 \mcd 3^2 \mcd 5 \mcd 7 \mcd 13$\\
$U_{4}(5)$ & $2^7 \mcd 3^4 \mcd 5^6 \mcd 7 \mcd 13$\\
$L_{3}(9)$ & $2^7 \mcd 3^6 \mcd 5 \mcd 7 \mcd 13$\\
$S_{6}(3)$ & $2^9 \mcd 3^9 \mcd 5 \mcd 7 \mcd 13$\\
$O_{7}(3)$ & $2^9 \mcd 3^9 \mcd 5 \mcd 7 \mcd 13$\\
$G_2(4)$ & $2^{12} \mcd 3^3 \mcd 5^2 \mcd 7 \mcd 13$\\
$S_{4}(8)$ & $2^{12} \mcd 3^4 \mcd 5 \mcd 7^2 \mcd 13$\\
$O^+_{8}(3)$ & $2^{12} \mcd 3^{12} \mcd 5^2 \mcd 7 \mcd 13$\\
$L_{5}(3)$ & $2^9 \mcd 3^{10} \mcd 5 \mcd 11^2 \mcd 13$\\
$A_{13}$ & $2^9 \mcd 3^5 \mcd 5^2 \mcd 7 \mcd 11 \mcd 13$\\
$A_{14}$ & $2^{10} \mcd 3^5 \mcd 5^2 \mcd 7^2 \mcd 11 \mcd 13$\\
$A_{15}$ & $2^{10} \mcd 3^6 \mcd 5^3 \mcd 7^2 \mcd 11 \mcd 13$\\
$L_{6}(3)$ & $2^{11} \mcd 3^{15} \mcd 5 \mcd 7 \mcd 11^2 \mcd 13^2$\\
$Suz$ & $2^{13} \mcd 3^7 \mcd 5^2 \mcd 7 \mcd 11 \mcd 13$\\
$A_{16}$ & $2^{14} \mcd 3^6 \mcd 5^3 \mcd 7^2 \mcd 11 \mcd 13$\\
$Fi_{22}$ & $2^{17} \mcd 3^9 \mcd 5^2 \mcd 7 \mcd 11 \mcd 13$\\
\hline
\multicolumn{2}{c}{$\S_{17}$ \hfill $^{\vphantom{A^A}}17\in\pi(G)\se\{2,3,5,\ld,17\}$ \hfill $\#\,18$}  \\
\hline
 $^{\vphantom{A^A}}$$L_{2}(17)$ & $2^4 \mcd 3^2 \mcd 17$\\
$L_{2}(16)$ & $2^4 \mcd 3 \mcd 5 \mcd 17$\\
$S_{4}(4)$ & $2^8 \mcd 3^2 \mcd 5^2 \mcd 17$\\
$He$ & $2^{10} \mcd 3^3 \mcd 5^2 \mcd 7^3 \mcd 17$\\
$O^-_{8}(2)$ & $2^{12} \mcd 3^4 \mcd 5 \mcd 7 \mcd 17$\\
$L_{4}(4)$ & $2^{12} \mcd 3^4 \mcd 5^2 \mcd 7 \mcd 17$\\
$S_{8}(2)$ & $2^{16} \mcd 3^5 \mcd 5^2 \mcd 7 \mcd 17$\\
$U_{4}(4)$ & $2^{12} \mcd 3^2 \mcd 5^3 \mcd 13 \mcd 17$\\
$U_{3}(17)$ & $2^6 \mcd 3^4 \mcd 7 \mcd 13 \mcd 17^3$\\
$O^-_{10}(2)$ & $2^{20} \mcd 3^6 \mcd 5^2 \mcd 7 \mcd 11 \mcd 17$\\
$L_{2}(13^{2})$ & $2^3 \mcd 3 \mcd 5 \mcd 7 \mcd 13^2 \mcd 17$\\
$S_{4}(13)$ & $2^6 \mcd 3^2 \mcd 5 \mcd 7^2 \mcd 13^4 \mcd 17$\\
$L_{3}(16)$ & $2^{12} \mcd 3^2 \mcd 5^2 \mcd 7 \mcd 13 \mcd 17$\\
$S_{6}(4)$ & $2^{18} \mcd 3^4 \mcd 5^3 \mcd 7 \mcd 13 \mcd 17$\\
$O^+_{8}(4)$ & $2^{24} \mcd 3^5 \mcd 5^4 \mcd 7 \mcd 13 \mcd 17^2$\\
$F_4(2)$ & $2^{24} \mcd 3^6 \mcd 5^2 \mcd 7^2 \mcd 13 \mcd 17$\\
$A_{17}$ & $2^{14} \mcd 3^6 \mcd 5^3 \mcd 7^2 \mcd 11 \mcd 13 \mcd 17$\\
$A_{18}$ & $2^{15} \mcd 3^8 \mcd 5^3 \mcd 7^2 \mcd 11 \mcd 13 \mcd 17$\\
\hline
\multicolumn{2}{c}{$\S_{19}$ \hfill $^{\vphantom{A^A}}19\in\pi(G)\se\{2,3,5,\ld,19\}$ \hfill $\#\,15$}  \\
\hline
 $^{\vphantom{A^A}}$$L_{2}(19)$ & $2^2 \mcd 3^2 \mcd 5 \mcd 19$\\
$L_{3}(7)$ & $2^5 \mcd 3^2 \mcd 7^3 \mcd 19$\\
$U_{3}(8)$ & $2^9 \mcd 3^4 \mcd 7 \mcd 19$\\
$U_{3}(19)$ & $2^5 \mcd 3^2 \mcd 5^2 \mcd 7^3 \mcd 19^3$\\
$L_{4}(7)$ & $2^9 \mcd 3^4 \mcd 5^2 \mcd 7^6 \mcd 19$\\
$J_3$ & $2^7 \mcd 3^5 \mcd 5 \mcd 17 \mcd 19$\\
$J_1$ & $2^3 \mcd 3 \mcd 5 \mcd 7 \mcd 11 \mcd 19$\\
$L_{3}(11)$ & $2^4 \mcd 3 \mcd 5^2 \mcd 7 \mcd 11^3 \mcd 19$\\
$HN$ & $2^{14} \mcd 3^6 \mcd 5^6 \mcd 7 \mcd 11 \mcd 19$\\
$U_{4}(8)$ & $2^{18} \mcd 3^7 \mcd 5 \mcd 7^2 \mcd 13 \mcd 19$\\
$A_{19}$ & $2^{15} \mcd 3^8 \mcd 5^3 \mcd 7^2 \mcd 11 \mcd 13 \mcd 17 \mcd 19$\\
$A_{20}$ & $2^{17} \mcd 3^8 \mcd 5^4 \mcd 7^2 \mcd 11 \mcd 13 \mcd 17 \mcd 19$\\
$A_{21}$ & $2^{17} \mcd 3^9 \mcd 5^4 \mcd 7^3 \mcd 11 \mcd 13 \mcd 17 \mcd 19$\\
$A_{22}$ & $2^{18} \mcd 3^9 \mcd 5^4 \mcd 7^3 \mcd 11^2 \mcd 13 \mcd 17 \mcd 19$\\
$^2E_6(2)$ & $2^{36} \mcd 3^9 \mcd 5^2 \mcd 7^2 \mcd 11 \mcd 13 \mcd 17 \mcd 19$\\
\hline
\multicolumn{2}{c}{$\S_{23}$ \hfill $^{\vphantom{A^A}}23\in\pi(G)\se\{2,3,5,\ld,23\}$ \hfill $\#\,14$}  \\
\hline
 $^{\vphantom{A^A}}$$L_{2}(23)$ & $2^3 \mcd 3 \mcd 11 \mcd 23$\\
$U_{3}(23)$ & $2^7 \mcd 3^2 \mcd 11 \mcd 13^2 \mcd 23^3$\\
$M_{23}$ & $2^7 \mcd 3^2 \mcd 5 \mcd 7 \mcd 11 \mcd 23$\\
$M_{24}$ & $2^{10} \mcd 3^3 \mcd 5 \mcd 7 \mcd 11 \mcd 23$\\
$Co_3$ & $2^{10} \mcd 3^7 \mcd 5^3 \mcd 7 \mcd 11 \mcd 23$\\
$Co_2$ & $2^{18} \mcd 3^6 \mcd 5^3 \mcd 7 \mcd 11 \mcd 23$\\
$Co_1$ & $2^{21} \mcd 3^9 \mcd 5^4 \mcd 7^2 \mcd 11 \mcd 13 \mcd 23$\\
$Fi_{23}$ & $2^{18} \mcd 3^{13} \mcd 5^2 \mcd 7 \mcd 11 \mcd 13 \mcd 17 \mcd 23$\\
$A_{23},\ld,A_{28}$ & $\pi(G)=\{2,3,5,\ld,23\}$\\
\hline
\multicolumn{2}{c}{$\S_{29}$ \hfill $^{\vphantom{A^A}}29\in\pi(G)\se\{2,3,5,\ld,29\}$ \hfill $\#\,8$}  \\
\hline
 $^{\vphantom{A^A}}$$L_{2}(29)$ & $2^2 \mcd 3 \mcd 5 \mcd 7 \mcd 29$\\
$L_{2}(17^{2})$ & $2^5 \mcd 3^2 \mcd 5 \mcd 17^2 \mcd 29$\\
$S_{4}(17)$ & $2^{10} \mcd 3^4 \mcd 5 \mcd 17^4 \mcd 29$\\
$Ru$ & $2^{14} \mcd 3^3 \mcd 5^3 \mcd 7 \mcd 13 \mcd 29$\\
$U_{4}(17)$ & $2^{11} \mcd 3^7 \mcd 5 \mcd 7 \mcd 13 \mcd 17^6 \mcd 29$\\
$Fi_{24}'$ & $2^{21} \mcd 3^{16} \mcd 5^2 \mcd 7^3 \mcd 11 \mcd 13 \mcd 17 \mcd 23 \mcd 29$\\
$A_{29},A_{30}$ & $\pi(G)=\{2,3,5,\ld,29\}$\\
\hline
\multicolumn{2}{c}{$\S_{31}$ \hfill $^{\vphantom{A^A}}31\in\pi(G)\se\{2,3,5,\ld,31\}$ \hfill $\#\,28$}  \\
\hline
 $^{\vphantom{A^A}}$$L_{2}(31)$ & $2^5 \mcd 3 \mcd 5 \mcd 31$\\
$L_{3}(5)$ & $2^5 \mcd 3 \mcd 5^3 \mcd 31$\\
$L_{2}(32)$ & $2^5 \mcd 3 \mcd 11 \mcd 31$\\
$L_{2}(5^{3})$ & $2^2 \mcd 3^2 \mcd 5^3 \mcd 7 \mcd 31$\\
$G_2(5)$ & $2^6 \mcd 3^3 \mcd 5^6 \mcd 7 \mcd 31$\\
$L_{5}(2)$ & $2^{10} \mcd 3^2 \mcd 5 \mcd 7 \mcd 31$\\
$L_{6}(2)$ & $2^{15} \mcd 3^4 \mcd 5 \mcd 7^2 \mcd 31$\\
$L_{4}(5)$ & $2^7 \mcd 3^2 \mcd 5^6 \mcd 13 \mcd 31$\\
$L_{3}(25)$ & $2^7 \mcd 3^2 \mcd 5^6 \mcd 7 \mcd 13 \mcd 31$\\
$O_{7}(5)$ & $2^9 \mcd 3^4 \mcd 5^9 \mcd 7 \mcd 13 \mcd 31$\\
$S_{6}(5)$ & $2^9 \mcd 3^4 \mcd 5^9 \mcd 7 \mcd 13 \mcd 31$\\
$O^+_{8}(5)$ & $2^{12} \mcd 3^5 \mcd 5^{12} \mcd 7 \mcd 13^2 \mcd 31$\\
$O^+_{10}(2)$ & $2^{20} \mcd 3^5 \mcd 5^2 \mcd 7 \mcd 17 \mcd 31$\\
$U_{3}(31)$ & $2^{11} \mcd 3 \mcd 5 \mcd 7^2 \mcd 19 \mcd 31^3$\\
$L_{5}(4)$ & $2^{20} \mcd 3^5 \mcd 5^2 \mcd 7 \mcd 11 \mcd 17 \mcd 31$\\
$S_{10}(2)$ & $2^{25} \mcd 3^6 \mcd 5^2 \mcd 7 \mcd 11 \mcd 17 \mcd 31$\\
$O^+_{12}(2)$ & $2^{30} \mcd 3^8 \mcd 5^2 \mcd 7^2 \mcd 11 \mcd 17 \mcd 31$\\
$ON$ & $2^9 \mcd 3^4 \mcd 5 \mcd 7^3 \mcd 11 \mcd 19 \mcd 31$\\
$Th$ & $2^{15} \mcd 3^{10} \mcd 5^3 \mcd 7^2 \mcd 13 \mcd 19 \mcd 31$\\
$O^-_{12}(2)$ & $2^{30} \mcd 3^6 \mcd 5^3 \mcd 7 \mcd 11 \mcd 13 \mcd 17 \mcd 31$\\
$L_{6}(4)$ & $2^{30} \mcd 3^6 \mcd 5^3 \mcd 7^2 \mcd 11 \mcd 13 \mcd 17 \mcd 31$\\
$S_{12}(2)$ & $2^{36} \mcd 3^8 \mcd 5^3 \mcd 7^2 \mcd 11 \mcd 13 \mcd 17 \mcd 31$\\
$A_{31},\ld,A_{36}$ & $\pi(G)=\{2,3,5,\ld,31\}$\\
\hline
\multicolumn{2}{c}{$\S_{37}$ \hfill $^{\vphantom{A^A}}37\in\pi(G)\se\{2,3,5,\ld,37\}$ \hfill $\#\,13$}  \\
\hline
 $^{\vphantom{A^A}}$$L_{2}(37)$ & $2^2 \mcd 3^2 \mcd 19 \mcd 37$\\
$U_{3}(11)$ & $2^5 \mcd 3^2 \mcd 5 \mcd 11^3 \mcd 37$\\
$L_{2}(31^{2})$ & $2^6 \mcd 3 \mcd 5 \mcd 13 \mcd 31^2 \mcd 37$\\
$S_{4}(31)$ & $2^{12} \mcd 3^2 \mcd 5^2 \mcd 13 \mcd 31^4 \mcd 37$\\
$^2G_2(27)$ & $2^3 \mcd 3^9 \mcd 7 \mcd 13 \mcd 19 \mcd 37$\\
$U_{3}(27)$ & $2^5 \mcd 3^9 \mcd 7^2 \mcd 13 \mcd 19 \mcd 37$\\
$L_{2}(11^{3})$ & $2^2 \mcd 3^2 \mcd 5 \mcd 7 \mcd 11^3 \mcd 19 \mcd 37$\\
$G_2(11)$ & $2^6 \mcd 3^3 \mcd 5^2 \mcd 7 \mcd 11^6 \mcd 19 \mcd 37$\\
$U_{4}(31)$ & $2^{16} \mcd 3^2 \mcd 5^2 \mcd 7^2 \mcd 13 \mcd 19 \mcd 31^6 \mcd 37$\\
$A_{37},\ld,A_{40}$ & $\pi(G)=\{2,3,5,\ld,37\}$\\
\hline
\multicolumn{2}{c}{$\S_{41}$ \hfill $^{\vphantom{A^A}}41\in\pi(G)\se\{2,3,5,\ld,41\}$ \hfill $\#\,17$}  \\
\hline
 $^{\vphantom{A^A}}$$L_{2}(3^{4})$ & $2^4 \mcd 3^4 \mcd 5 \mcd 41$\\
$S_{4}(9)$ & $2^8 \mcd 3^8 \mcd 5^2 \mcd 41$\\
$Sz(32)$ & $2^{10} \mcd 5^2 \mcd 31 \mcd 41$\\
$L_{2}(41)$ & $2^3 \mcd 3 \mcd 5 \mcd 7 \mcd 41$\\
$O^-_{8}(3)$ & $2^{10} \mcd 3^{12} \mcd 5 \mcd 7 \mcd 13 \mcd 41$\\
$L_{4}(9)$ & $2^{10} \mcd 3^{12} \mcd 5^2 \mcd 7 \mcd 13 \mcd 41$\\
$S_{8}(3)$ & $2^{14} \mcd 3^{16} \mcd 5^2 \mcd 7 \mcd 13 \mcd 41$\\
$O_{9}(3)$ & $2^{14} \mcd 3^{16} \mcd 5^2 \mcd 7 \mcd 13 \mcd 41$\\
$L_{2}(41^{2})$ & $2^4 \mcd 3 \mcd 5 \mcd 7 \mcd 29^2 \mcd 41^2$\\
$S_{4}(41)$ & $2^8 \mcd 3^2 \mcd 5^2 \mcd 7^2 \mcd 29^2 \mcd 41^4$\\
$L_{2}(2^{10})$ & $2^{10} \mcd 3 \mcd 5^2 \mcd 11 \mcd 31 \mcd 41$\\
$S_{4}(32)$ & $2^{20} \mcd 3^2 \mcd 5^2 \mcd 11^2 \mcd 31^2 \mcd 41$\\
$U_{5}(4)$ & $2^{20} \mcd 3^2 \mcd 5^4 \mcd 13 \mcd 17 \mcd 41$\\
$O^+_{10}(3)$ & $2^{15} \mcd 3^{20} \mcd 5^2 \mcd 7 \mcd 11^2 \mcd 13 \mcd 41$\\
$U_{6}(4)$ & $2^{30} \mcd 3^4 \mcd 5^6 \mcd 7 \mcd 13^2 \mcd 17 \mcd 41$\\
$A_{41},A_{42}$ & $\pi(G)=\{2,3,5,\ld,41\}$\\
\hline
\multicolumn{2}{c}{$\S_{43}$ \hfill $^{\vphantom{A^A}}43\in\pi(G)\se\{2,3,5,\ld,43\}$ \hfill $\#\,22$}  \\
\hline
 $^{\vphantom{A^A}}$$U_{3}(7)$ & $2^7 \mcd 3 \mcd 7^3 \mcd 43$\\
$U_{4}(7)$ & $2^{10} \mcd 3^2 \mcd 5^2 \mcd 7^6 \mcd 43$\\
$L_{2}(43)$ & $2^2 \mcd 3 \mcd 7 \mcd 11 \mcd 43$\\
$L_{2}(7^{3})$ & $2^3 \mcd 3^2 \mcd 7^3 \mcd 19 \mcd 43$\\
$G_2(7)$ & $2^8 \mcd 3^3 \mcd 7^6 \mcd 19 \mcd 43$\\
$U_{7}(2)$ & $2^{21} \mcd 3^8 \mcd 5 \mcd 7 \mcd 11 \mcd 43$\\
$L_{3}(49)$ & $2^9 \mcd 3^2 \mcd 5^2 \mcd 7^6 \mcd 19 \mcd 43$\\
$S_{6}(7)$ & $2^{12} \mcd 3^4 \mcd 5^2 \mcd 7^9 \mcd 19 \mcd 43$\\
$O_{7}(7)$ & $2^{12} \mcd 3^4 \mcd 5^2 \mcd 7^9 \mcd 19 \mcd 43$\\
$O^+_{8}(7)$ & $2^{16} \mcd 3^5 \mcd 5^4 \mcd 7^{12} \mcd 19 \mcd 43$\\
$U_{3}(37)$ & $2^4 \mcd 3^2 \mcd 19^2 \mcd 31 \mcd 37^3 \mcd 43$\\
$U_{8}(2)$ & $2^{28} \mcd 3^9 \mcd 5^2 \mcd 7 \mcd 11 \mcd 17 \mcd 43$\\
$L_{2}(43^{2})$ & $2^3 \mcd 3 \mcd 5^2 \mcd 7 \mcd 11 \mcd 37 \mcd 43^2$\\
$S_{4}(43)$ & $2^6 \mcd 3^2 \mcd 5^2 \mcd 7^2 \mcd 11^2 \mcd 37 \mcd 43^4$\\
$U_{9}(2)$ & $2^{36} \mcd 3^{11} \mcd 5^2 \mcd 7 \mcd 11 \mcd 17 \mcd 19 \mcd 43$\\
$O^-_{14}(2)$ & $2^{42} \mcd 3^9 \mcd 5^3 \mcd 7^2 \mcd 11 \mcd 13 \mcd 17 \mcd 31 \mcd 43$\\
$U_{10}(2)$ & $2^{45} \mcd 3^{13} \mcd 5^2 \mcd 7 \mcd 11^2 \mcd 17 \mcd 19 \mcd 31 \mcd 43$\\
$J_4$ & $2^{21} \mcd 3^3 \mcd 5 \mcd 7 \mcd 11^3 \mcd 23 \mcd 29 \mcd 31 \mcd 37 \mcd 43$\\
$A_{43},\ld,A_{46}$ & $\pi(G)=\{2,3,5,\ld,43\}$\\
\hline
\multicolumn{2}{c}{$\S_{47}$ \hfill $^{\vphantom{A^A}}47\in\pi(G)\se\{2,3,5,\ld,47\}$ \hfill $\#\,10$}  \\
\hline
 $^{\vphantom{A^A}}$$L_{2}(47)$ & $2^4 \mcd 3 \mcd 23 \mcd 47$\\
$L_{2}(47^{2})$ & $2^5 \mcd 3 \mcd 5 \mcd 13 \mcd 17 \mcd 23 \mcd 47^2$\\
$S_{4}(47)$ & $2^{10} \mcd 3^2 \mcd 5 \mcd 13 \mcd 17 \mcd 23^2 \mcd 47^4$\\
$B$ & $2^{41} \mcd 3^{13} \mcd 5^6 \mcd 7^2 \mcd 11 \mcd 13 \mcd 17 \mcd 19 \mcd 23 \mcd 31 \mcd 47$\\
$A_{47},\ld,A_{52}$ & $\pi(G)=\{2,3,5,\ld,47\}$\\
\hline
\multicolumn{2}{c}{$\S_{53}$ \hfill $^{\vphantom{A^A}}53\in\pi(G)\se\{2,3,5,\ld,53\}$ \hfill $\#\,10$}  \\
\hline
 $^{\vphantom{A^A}}$$L_{2}(53)$ & $2^2 \mcd 3^3 \mcd 13 \mcd 53$\\
$L_{2}(23^{2})$ & $2^4 \mcd 3 \mcd 5 \mcd 11 \mcd 23^2 \mcd 53$\\
$S_{4}(23)$ & $2^8 \mcd 3^2 \mcd 5 \mcd 11^2 \mcd 23^4 \mcd 53$\\
$U_{4}(23)$ & $2^{10} \mcd 3^4 \mcd 5 \mcd 11^2 \mcd 13^2 \mcd 23^6 \mcd 53$\\
$A_{53},\ld,A_{58}$ & $\pi(G)=\{2,3,5,\ld,53\}$\\
\hline
\multicolumn{2}{c}{$\S_{59}$ \hfill $^{\vphantom{A^A}}59\in\pi(G)\se\{2,3,5,\ld,59\}$ \hfill $\#\,3$}  \\
\hline
 $^{\vphantom{A^A}}$$L_{2}(59)$ & $2^2 \mcd 3 \mcd 5 \mcd 29 \mcd 59$\\
$A_{59},A_{60}$ & $\pi(G)=\{2,3,5,\ld,59\}$\\
\hline
\multicolumn{2}{c}{$\S_{61}$ \hfill $^{\vphantom{A^A}}61\in\pi(G)\se\{2,3,5,\ld,61\}$ \hfill $\#\,27$}  \\
\hline
 $^{\vphantom{A^A}}$$L_{2}(3^{5})$ & $2^2 \mcd 3^5 \mcd 11^2 \mcd 61$\\
$U_{5}(3)$ & $2^{11} \mcd 3^{10} \mcd 5 \mcd 7 \mcd 61$\\
$L_{2}(11^{2})$ & $2^3 \mcd 3 \mcd 5 \mcd 11^2 \mcd 61$\\
$S_{4}(11)$ & $2^6 \mcd 3^2 \mcd 5^2 \mcd 11^4 \mcd 61$\\
$L_{2}(61)$ & $2^2 \mcd 3 \mcd 5 \mcd 31 \mcd 61$\\
$L_{3}(13)$ & $2^5 \mcd 3^2 \mcd 7 \mcd 13^3 \mcd 61$\\
$U_{6}(3)$ & $2^{13} \mcd 3^{15} \mcd 5 \mcd 7^2 \mcd 13 \mcd 61$\\
$U_{4}(11)$ & $2^7 \mcd 3^4 \mcd 5^2 \mcd 11^6 \mcd 37 \mcd 61$\\
$L_{3}(47)$ & $2^6 \mcd 3 \mcd 23^2 \mcd 37 \mcd 47^3 \mcd 61$\\
$L_{4}(11)$ & $2^7 \mcd 3^2 \mcd 5^3 \mcd 7 \mcd 11^6 \mcd 19 \mcd 61$\\
$L_{4}(13)$ & $2^7 \mcd 3^4 \mcd 5 \mcd 7^2 \mcd 13^6 \mcd 17 \mcd 61$\\
$O^-_{10}(3)$ & $2^{15} \mcd 3^{20} \mcd 5^2 \mcd 7 \mcd 13 \mcd 41 \mcd 61$\\
$L_{5}(9)$ & $2^{15} \mcd 3^{20} \mcd 5^2 \mcd 7 \mcd 11^2 \mcd 13 \mcd 41 \mcd 61$\\
$S_{10}(3)$ & $2^{17} \mcd 3^{25} \mcd 5^2 \mcd 7 \mcd 11^2 \mcd 13 \mcd 41 \mcd 61$\\
$O_{11}(3)$ & $2^{17} \mcd 3^{25} \mcd 5^2 \mcd 7 \mcd 11^2 \mcd 13 \mcd 41 \mcd 61$\\
$O^+_{12}(3)$ & $2^{19} \mcd 3^{30} \mcd 5^2 \mcd 7^2 \mcd 11^2 \mcd 13^2 \mcd 41 \mcd 61$\\
$L_{3}(11^{2})$ & $2^7 \mcd 3^2 \mcd 5^2 \mcd 7 \mcd 11^6 \mcd 19 \mcd 37 \mcd 61$\\
$S_{6}(11)$ & $2^9 \mcd 3^4 \mcd 5^3 \mcd 7 \mcd 11^9 \mcd 19 \mcd 37 \mcd 61$\\
$O_{7}(11)$ & $2^9 \mcd 3^4 \mcd 5^3 \mcd 7 \mcd 11^9 \mcd 19 \mcd 37 \mcd 61$\\
$O^+_{8}(11)$ & $2^{12} \mcd 3^5 \mcd 5^4 \mcd 7 \mcd 11^{12} \mcd 19 \mcd 37 \mcd 61^2$\\
$L_{4}(47)$ & $2^{11} \mcd 3^2 \mcd 5 \mcd 13 \mcd 17 \mcd 23^3 \mcd 37 \mcd 47^6 \mcd 61$\\
$A_{61},\ld,A_{66}$ & $\pi(G)=\{2,3,5,\ld,61\}$\\
\hline
\multicolumn{2}{c}{$\S_{67}$ \hfill $^{\vphantom{A^A}}67\in\pi(G)\se\{2,3,5,\ld,67\}$ \hfill $\#\,11$}  \\
\hline
 $^{\vphantom{A^A}}$$L_{2}(67)$ & $2^2 \mcd 3 \mcd 11 \mcd 17 \mcd 67$\\
$L_{3}(37)$ & $2^5 \mcd 3^4 \mcd 7 \mcd 19 \mcd 37^3 \mcd 67$\\
$L_{3}(29)$ & $2^5 \mcd 3 \mcd 5 \mcd 7^2 \mcd 13 \mcd 29^3 \mcd 67$\\
$L_{3}(67)$ & $2^4 \mcd 3^2 \mcd 7^2 \mcd 11^2 \mcd 17 \mcd 31 \mcd 67^3$\\
$Ly$ & $2^8 \mcd 3^7 \mcd 5^6 \mcd 7 \mcd 11 \mcd 31 \mcd 37 \mcd 67$\\
$L_{2}(37^{3})$ & $2^2 \mcd 3^3 \mcd 7 \mcd 19 \mcd 31 \mcd 37^3 \mcd 43 \mcd 67$\\
$G_2(37)$ & $2^6 \mcd 3^5 \mcd 7 \mcd 19^2 \mcd 31 \mcd 37^6 \mcd 43 \mcd 67$\\
$A_{67},\ld,A_{70}$ & $\pi(G)=\{2,3,5,\ld,67\}$\\
\hline
\multicolumn{2}{c}{$\S_{71}$ \hfill $^{\vphantom{A^A}}71\in\pi(G)\se\{2,3,5,\ld,71\}$ \hfill $\#\,6$}  \\
\hline
 $^{\vphantom{A^A}}$$L_{2}(71)$ & $2^3 \mcd 3^2 \mcd 5 \mcd 7 \mcd 71$\\
$L_{5}(5)$ & $2^{11} \mcd 3^2 \mcd 5^{10} \mcd 11 \mcd 13 \mcd 31 \mcd 71$\\
$L_{6}(5)$ & $2^{13} \mcd 3^4 \mcd 5^{15} \mcd 7 \mcd 11 \mcd 13 \mcd 31^2 \mcd 71$\\
$M$ & $2^{46} \mcd 3^{20} \mcd 5^9 \mcd 7^6 \mcd 11^2 \mcd 13^3 \mcd 17 \mcd 19 \mcd 23 \mcd 29 \mcd 31 \mcd 41 \mcd 47 \mcd 59 \mcd 71$\\
$A_{71},A_{72}$ & $\pi(G)=\{2,3,5,\ld,71\}$\\
\hline
\multicolumn{2}{c}{$\S_{73}$ \hfill $^{\vphantom{A^A}}73\in\pi(G)\se\{2,3,5,\ld,73\}$ \hfill $\#\,34$}  \\
\hline
 $^{\vphantom{A^A}}$$U_{3}(9)$ & $2^5 \mcd 3^6 \mcd 5^2 \mcd 73$\\
$L_{3}(8)$ & $2^9 \mcd 3^2 \mcd 7^2 \mcd 73$\\
$L_{2}(73)$ & $2^3 \mcd 3^2 \mcd 37 \mcd 73$\\
$U_{4}(9)$ & $2^9 \mcd 3^{12} \mcd 5^3 \mcd 41 \mcd 73$\\
$^3D_4(3)$ & $2^6 \mcd 3^{12} \mcd 7^2 \mcd 13^2 \mcd 73$\\
$L_{2}(2^{9})$ & $2^9 \mcd 3^3 \mcd 7 \mcd 19 \mcd 73$\\
$G_2(8)$ & $2^{18} \mcd 3^5 \mcd 7^2 \mcd 19 \mcd 73$\\
$L_{2}(3^{6})$ & $2^3 \mcd 3^6 \mcd 5 \mcd 7 \mcd 13 \mcd 73$\\
$S_{4}(27)$ & $2^6 \mcd 3^{12} \mcd 5 \mcd 7^2 \mcd 13^2 \mcd 73$\\
$G_2(9)$ & $2^8 \mcd 3^{12} \mcd 5^2 \mcd 7 \mcd 13 \mcd 73$\\
$L_{4}(8)$ & $2^{18} \mcd 3^4 \mcd 5 \mcd 7^3 \mcd 13 \mcd 73$\\
$L_{3}(64)$ & $2^{18} \mcd 3^4 \mcd 5 \mcd 7^2 \mcd 13 \mcd 19 \mcd 73$\\
$S_{6}(8)$ & $2^{27} \mcd 3^7 \mcd 5 \mcd 7^3 \mcd 13 \mcd 19 \mcd 73$\\
$O^+_{8}(8)$ & $2^{36} \mcd 3^9 \mcd 5^2 \mcd 7^4 \mcd 13^2 \mcd 19 \mcd 73$\\
$L_{3}(3^{4})$ & $2^9 \mcd 3^{12} \mcd 5^2 \mcd 7 \mcd 13 \mcd 41 \mcd 73$\\
$S_{6}(9)$ & $2^{12} \mcd 3^{18} \mcd 5^3 \mcd 7 \mcd 13 \mcd 41 \mcd 73$\\
$O_{7}(9)$ & $2^{12} \mcd 3^{18} \mcd 5^3 \mcd 7 \mcd 13 \mcd 41 \mcd 73$\\
$F_4(3)$ & $2^{15} \mcd 3^{24} \mcd 5^2 \mcd 7^2 \mcd 13^2 \mcd 41 \mcd 73$\\
$O^+_{8}(9)$ & $2^{16} \mcd 3^{24} \mcd 5^4 \mcd 7 \mcd 13 \mcd 41^2 \mcd 73$\\
$L_{2}(73^{2})$ & $2^4 \mcd 3^2 \mcd 5 \mcd 13 \mcd 37 \mcd 41 \mcd 73^2$\\
$S_{4}(73)$ & $2^8 \mcd 3^4 \mcd 5 \mcd 13 \mcd 37^2 \mcd 41 \mcd 73^4$\\
$E_6(2)$ & $2^{36} \mcd 3^6 \mcd 5^2 \mcd 7^3 \mcd 13 \mcd 17 \mcd 31 \mcd 73$\\
$U_{4}(27)$ & $2^7 \mcd 3^{18} \mcd 5 \mcd 7^3 \mcd 13^2 \mcd 19 \mcd 37 \mcd 73$\\
$O^-_{12}(3)$ & $2^{18} \mcd 3^{30} \mcd 5^3 \mcd 7 \mcd 11^2 \mcd 13 \mcd 41 \mcd 61 \mcd 73$\\
$L_{6}(9)$ & $2^{18} \mcd 3^{30} \mcd 5^3 \mcd 7^2 \mcd 11^2 \mcd 13^2 \mcd 41 \mcd 61 \mcd 73$\\
$O_{13}(3)$ & $2^{21} \mcd 3^{36} \mcd 5^3 \mcd 7^2 \mcd 11^2 \mcd 13^2 \mcd 41 \mcd 61 \mcd 73$\\
$S_{12}(3)$ & $2^{21} \mcd 3^{36} \mcd 5^3 \mcd 7^2 \mcd 11^2 \mcd 13^2 \mcd 41 \mcd 61 \mcd 73$\\
$^2E_6(3)$ & $2^{19} \mcd 3^{36} \mcd 5^2 \mcd 7^3 \mcd 13^2 \mcd 19 \mcd 37 \mcd 41 \mcd 61 \mcd 73$\\
$A_{73},\ld,A_{78}$ & $\pi(G)=\{2,3,5,\ld,73\}$\\
\hline
\multicolumn{2}{c}{$\S_{79}$ \hfill $^{\vphantom{A^A}}79\in\pi(G)\se\{2,3,5,\ld,79\}$ \hfill $\#\,14$}  \\
\hline
 $^{\vphantom{A^A}}$$L_{2}(79)$ & $2^4 \mcd 3 \mcd 5 \mcd 13 \mcd 79$\\
$L_{3}(23)$ & $2^5 \mcd 3 \mcd 7 \mcd 11^2 \mcd 23^3 \mcd 79$\\
$L_{3}(79)$ & $2^6 \mcd 3^2 \mcd 5 \mcd 7^2 \mcd 13^2 \mcd 43 \mcd 79^3$\\
$L_{2}(23^{3})$ & $2^3 \mcd 3^2 \mcd 7 \mcd 11 \mcd 13^2 \mcd 23^3 \mcd 79$\\
$G_2(23)$ & $2^8 \mcd 3^3 \mcd 7 \mcd 11^2 \mcd 13^2 \mcd 23^6 \mcd 79$\\
$L_{4}(23)$ & $2^9 \mcd 3^2 \mcd 5 \mcd 7 \mcd 11^3 \mcd 23^6 \mcd 53 \mcd 79$\\
$L_{3}(23^{2})$ & $2^9 \mcd 3^2 \mcd 5 \mcd 7 \mcd 11^2 \mcd 13^2 \mcd 23^6 \mcd 53 \mcd 79$\\
$S_{6}(23)$ & $2^{12} \mcd 3^4 \mcd 5 \mcd 7 \mcd 11^3 \mcd 13^2 \mcd 23^9 \mcd 53 \mcd 79$\\
$O_{7}(23)$ & $2^{12} \mcd 3^4 \mcd 5 \mcd 7 \mcd 11^3 \mcd 13^2 \mcd 23^9 \mcd 53 \mcd 79$\\
$O^+_{8}(23)$ & $2^{16} \mcd 3^5 \mcd 5^2 \mcd 7 \mcd 11^4 \mcd 13^2 \mcd 23^{12} \mcd 53^2 \mcd 79$\\
$A_{79},\ld,A_{82}$ & $\pi(G)=\{2,3,5,\ld,79\}$\\
\hline
\multicolumn{2}{c}{$\S_{83}$ \hfill $^{\vphantom{A^A}}83\in\pi(G)\se\{2,3,5,\ld,83\}$ \hfill $\#\,9$}  \\
\hline
 $^{\vphantom{A^A}}$$L_{2}(83)$ & $2^2 \mcd 3 \mcd 7 \mcd 41 \mcd 83$\\
$L_{2}(83^{2})$ & $2^3 \mcd 3 \mcd 5 \mcd 7 \mcd 13 \mcd 41 \mcd 53 \mcd 83^2$\\
$S_{4}(83)$ & $2^6 \mcd 3^2 \mcd 5 \mcd 7^2 \mcd 13 \mcd 41^2 \mcd 53 \mcd 83^4$\\
$A_{83},\ld,A_{88}$ & $\pi(G)=\{2,3,5,\ld,83\}$\\
\hline
\multicolumn{2}{c}{$\S_{89}$ \hfill $^{\vphantom{A^A}}89\in\pi(G)\se\{2,3,5,\ld,89\}$ \hfill $\#\,9$}  \\
\hline
 $^{\vphantom{A^A}}$$L_{2}(89)$ & $2^3 \mcd 3^2 \mcd 5 \mcd 11 \mcd 89$\\
$A_{89},\ld,A_{96}$ & $\pi(G)=\{2,3,5,\ld,89\}$\\
\hline
\multicolumn{2}{c}{$\S_{97}$ \hfill $^{\vphantom{A^A}}97\in\pi(G)\se\{2,3,5,\ld,97\}$ \hfill $\#\,6$}  \\
\hline
 $^{\vphantom{A^A}}$$L_{2}(97)$ & $2^5 \mcd 3 \mcd 7^2 \mcd 97$\\
$L_{3}(61)$ & $2^5 \mcd 3^2 \mcd 5^2 \mcd 13 \mcd 31 \mcd 61^3 \mcd 97$\\
$A_{97},\ld,A_{100}$ & $\pi(G)=\{2,3,5,\ld,97\}$\\
\hline
\end{longtable}%
}

\bigskip
\medskip

\footnotesize{
\begin{longtable}[htb]{c}
\caption{\bf Primes $p\in\{100,\ld,1000\}$ with generic $\S_p$ }\label{gen}\\
\hline
$^{\vphantom{A^A}}$107, 131, 167, 197, 223, 227, 311, 317, 347, 359, 379,\\
383, 389, 397, 419, 461, 479, 503, 541, 569, 587, 617, \\
643, 647, 691, 761, 797, 827, 839, 887, 967, 977, 983 \\
\hline
\end{longtable}
}

\bigskip

\footnotesize{
\begin{longtable}[h]{l|l|l}
\caption {\bf Non-generic groups $G$ with  $100<p\in\pi(G)\se\{2,3,5,\ld,p\}$}\label{ng}\\
\hline
p & $|\S_p|$ & \multicolumn{1}{c}{$^{\vphantom{A^A}}G$}\\
\hline
$^{\vphantom{A^A}}$$101$ & $5$ & $U_{3}(101)$,\,$U_{5}(17)$\\
$103$ & $14$ & $U_{3}(47)$,\,$U_{3}(103)$,\,$L_{2}(47^{3})$,\,$G_2(47)$,\,$U_{4}(47)$,\,$L_{3}(47^{2})$,\,$S_{6}(47)$,\,$O_{7}(47)$,\\
 & &$O^+_{8}(47)$\\
$109$ & $12$ & $U_{3}(64)$,\,$^3D_4(8)$,\,$Sz(2^{9})$,\,$^2F_4(8)$,\,$L_{2}(2^{18})$,\,$G_2(64)$,\,$S_{4}(2^{9})$\\
$113$ & $16$ & $U_{7}(4)$\\
$127$ & $24$ & $L_{2}(2^{7})$,\,$L_{3}(19)$,\,$Sz(2^{7})$,\,$L_{2}(19^{3})$,\,$G_2(19)$,\,$L_{7}(2)$,\,$L_{8}(2)$,\,$L_{2}(2^{14})$,\\
 & &$S_{4}(2^{7})$,\,$L_{3}(107)$,\,$L_{9}(2)$,\,$O^+_{14}(2)$,\,$L_{10}(2)$,\,$L_{7}(4)$,\,$S_{14}(2)$,\,$O^+_{16}(2)$,\\
 & &$L_{11}(2)$,\,$E_7(2)$,\,$L_{12}(2)$\\
$137$ & $12$ & $L_{2}(37^{2})$,\,$S_{4}(37)$,\,$L_{4}(37)$,\,$U_{4}(37)$,\,$L_{3}(137)$,\,$L_{3}(37^{2})$,\,$S_{6}(37)$,\,$O_{7}(37)$,\\
 & &$O^+_{8}(37)$\\
$139$ & $14$ & $U_{3}(97)$,\,$U_{3}(43)$,\,$U_{4}(43)$\\
$149$ & $4$ & $L_{3}(149)$\\
$151$ & $11$ & $L_{3}(32)$,\,$L_{4}(32)$,\,$L_{5}(8)$,\,$L_{6}(8)$\\
$157$ & $17$ & $U_{3}(13)$,\,$L_{2}(13^{3})$,\,$G_2(13)$,\,$U_{4}(13)$,\,$L_{3}(13^{2})$,\,$S_{6}(13)$,\,$O_{7}(13)$,\,$O^+_{8}(13)$,\\
 & &$L_{2}(157^{2})$,\,$S_{4}(157)$\\
$163$ & $7$ & $U_{3}(59)$,\,$L_{3}(163)$\\
$173$ & $11$ & $L_{2}(173^{2})$,\,$S_{4}(173)$,\,$U_{3}(173)$,\,$U_{4}(173)$\\
$179$ & $4$ & $U_{3}(179)$\\
$181$ & $25$ & $L_{2}(19^{2})$,\,$S_{4}(19)$,\,$U_{3}(49)$,\,$U_{4}(19)$,\,$L_{4}(19)$,\,$L_{3}(19^{2})$,\,$O_{7}(19)$,\,$S_{6}(19)$,\\
 & &$O^+_{8}(19)$,\,$^3D_4(7)$,\,$L_{2}(7^{6})$,\,$S_{4}(7^{3})$,\,$G_2(49)$,\,$L_{3}(181)$\\
$191$ & $9$ & $U_{5}(7)$,\,$U_{6}(7)$,\,$L_{3}(191)$,\,$L_{2}(191^{2})$,\,$S_{4}(191)$,\,$L_{4}(191)$\\
$193$ & $15$ & $L_{2}(3^{8})$,\,$S_{4}(3^{4})$,\,$L_{2}(193^{2})$,\,$S_{4}(193)$,\,$U_{3}(109)$,\,$O^-_{8}(9)$,\,$L_{4}(3^{4})$,\,$S_{8}(9)$,\\
 & &$O_{9}(9)$,\,$O^+_{10}(9)$\\
$199$ & $16$ & $U_{3}(107)$,\,$L_{2}(107^{3})$,\,$G_2(107)$\\
$211$ & $20$ & $U_{3}(197)$,\,$L_{2}(211^{2})$,\,$S_{4}(211)$,\,$L_{3}(211)$,\,$U_{5}(23)$,\,$L_{4}(211)$,\,$U_{6}(23)$\\
$229$ & $14$ & $L_{2}(107^{2})$,\,$S_{4}(107)$,\,$L_{3}(229)$,\,$U_{4}(107)$,\,$L_{4}(107)$,\,$L_{3}(107^{2})$,\,$S_{6}(107)$,\,$O_{7}(107)$,\\
 & &$O^+_{8}(107)$\\
$233$ & $11$ & $L_{2}(89^{2})$,\,$S_{4}(89)$,\,$L_{2}(233^{2})$,\,$S_{4}(233)$\\
$239$ & $5$ & $L_{2}(239^{2})$,\,$S_{4}(239)$\\
$241$ & $25$ & $U_{3}(16)$,\,$^3D_4(4)$,\,$L_{2}(2^{12})$,\,$G_2(16)$,\,$S_{4}(64)$,\,$O^-_{8}(8)$,\,$L_{4}(64)$,\,$S_{8}(8)$,\\
 & &$U_{4}(64)$,\,$O^+_{10}(8)$,\,$L_{3}(2^{12})$,\,$S_{6}(64)$,\,$O^+_{8}(64)$,\,$F_4(8)$\\
$251$ & $9$ & $L_{2}(251^{2})$,\,$S_{4}(251)$\\
$257$ & $51$ & $L_{2}(2^{8})$,\,$S_{4}(16)$,\,$U_{4}(16)$,\,$O^-_{8}(4)$,\,$L_{4}(16)$,\,$S_{8}(4)$,\,$L_{2}(241^{2})$,\,$S_{4}(241)$,\\
 & &$U_{3}(257)$,\,$O^-_{10}(4)$,\,$L_{3}(2^{8})$,\,$S_{6}(16)$,\,$O^+_{8}(16)$,\,$F_4(4)$,\,$O^+_{10}(4)$,\,$L_{5}(16)$,\\
 & &$S_{10}(4)$,\,$O^+_{12}(4)$,\,$U_{8}(4)$,\,$O^-_{12}(4)$,\,$L_{6}(16)$,\,$S_{12}(4)$,\,$O^-_{16}(2)$,\,$L_{8}(4)$,\\
 & &$S_{16}(2)$,\,$^2E_6(4)$,\,$O^-_{18}(2)$,\,$E_6(4)$,\,$O^+_{18}(2)$,\,$U_{9}(4)$,\,$L_{9}(4)$,\,$S_{18}(2)$,\\
 & &$O^+_{20}(2)$,\,$O^-_{14}(4)$,\,$O^+_{14}(4)$,\,$O^-_{20}(2)$,\,$L_{10}(4)$,\,$S_{20}(2)$,\,$U_{10}(4)$,\,$L_{7}(16)$,\\
 & &$S_{14}(4)$,\,$O^+_{16}(4)$,\,$O^+_{22}(2)$,\,$E_7(4)$\\
$263$ & $11$ & $U_{3}(263)$,\,$L_{3}(263)$,\,$L_{2}(263^{3})$,\,$G_2(263)$\\
$269$ & $4$ & $L_{3}(269)$\\
$271$ & $13$ & $U_{3}(29)$,\,$^2G_2(3^{5})$,\,$U_{3}(3^{5})$,\,$L_{2}(29^{3})$,\,$G_2(29)$,\,$U_{5}(27)$\\
$277$ & $6$ & $L_{3}(277)$\\
$281$ & $5$ & $L_{2}(53^{2})$,\,$S_{4}(53)$\\
$283$ & $13$ & $U_{3}(239)$,\,$U_{4}(239)$\\
$293$ & $19$ & $U_{3}(293)$,\,$L_{2}(293^{2})$,\,$S_{4}(293)$,\,$U_{4}(293)$\\
$307$ & $16$ & $L_{3}(17)$,\,$L_{4}(17)$,\,$L_{2}(17^{3})$,\,$G_2(17)$,\,$L_{3}(17^{2})$,\,$S_{6}(17)$,\,$O_{7}(17)$,\,$O^+_{8}(17)$,\\
 & &$L_{2}(307^{2})$,\,$S_{4}(307)$,\,$U_{6}(17)$\\
$313$ & $17$ & $L_{2}(5^{4})$,\,$S_{4}(25)$,\,$L_{3}(313)$,\,$O^-_{8}(5)$,\,$L_{4}(25)$,\,$S_{8}(5)$,\,$O_{9}(5)$,\,$L_{2}(313^{2})$,\\
 & &$S_{4}(313)$,\,$O^+_{10}(5)$,\,$L_{4}(313)$,\,$^3D_4(29)$\\
$331$ & $32$ & $L_{3}(31)$,\,$U_{3}(32)$,\,$L_{2}(31^{3})$,\,$G_2(31)$,\,$U_{4}(32)$,\,$L_{4}(31)$,\,$L_{2}(2^{15})$,\,$G_2(32)$,\\
 & &$U_{5}(8)$,\,$U_{6}(8)$,\,$L_{3}(2^{10})$,\,$S_{6}(32)$,\,$O^+_{8}(32)$,\,$L_{3}(31^{2})$,\,$O_{7}(31)$,\,$S_{6}(31)$,\\
 & &$O^+_{8}(31)$,\,$O^-_{10}(8)$,\,$L_{5}(64)$,\,$S_{10}(8)$,\,$O^+_{12}(8)$,\,$O^-_{12}(8)$,\,$L_{6}(64)$,\,$S_{12}(8)$,\\
 & &$E_8(2)$\\
$337$ & $18$ & $L_{3}(2^{7})$,\,$L_{2}(337^{2})$,\,$S_{4}(337)$,\,$L_{4}(2^{7})$,\,$L_{7}(8)$,\,$L_{8}(8)$,\,$O^+_{14}(8)$\\
$349$ & $6$ & $U_{3}(227)$\\
$353$ & $10$ & $L_{2}(311^{2})$,\,$S_{4}(311)$,\,$U_{3}(353)$\\
$367$ & $10$ & $L_{3}(83)$,\,$L_{3}(283)$,\,$L_{4}(83)$\\
$373$ & $16$ & $U_{3}(89)$,\,$L_{2}(269^{2})$,\,$S_{4}(269)$,\,$U_{4}(89)$,\,$L_{3}(373)$,\,$U_{3}(373)$,\,$L_{4}(269)$,\,$L_{2}(373^{3})$,\\
 & &$G_2(373)$\\
$401$ & $11$ & $L_{2}(401^{2})$,\,$S_{4}(401)$\\
$409$ & $14$ & $L_{3}(53)$,\,$L_{4}(53)$,\,$^3D_4(49)$\\
$421$ & $29$ & $L_{2}(29^{2})$,\,$S_{4}(29)$,\,$U_{4}(29)$,\,$U_{3}(401)$,\,$L_{4}(29)$,\,$L_{2}(421^{2})$,\,$S_{4}(421)$,\,$L_{3}(29^{2})$,\\
 & &$O_{7}(29)$,\,$S_{6}(29)$,\,$O^+_{8}(29)$,\,$U_{3}(29^{2})$,\,$U_{5}(29)$,\,$U_{4}(401)$,\,$U_{6}(29)$,\,$L_{2}(29^{6})$,\\
 & &$S_{4}(29^{3})$,\,$G_2(29^{2})$\\
$431$ & $15$ & $U_{3}(431)$,\,$L_{2}(431^{2})$,\,$S_{4}(431)$,\,$L_{3}(431)$,\,$U_{4}(431)$,\,$L_{2}(431^{3})$,\,$G_2(431)$,\,$L_{4}(431)$,\\
 & &$L_{3}(431^{2})$,\,$S_{6}(431)$,\,$O_{7}(431)$,\,$O^+_{8}(431)$\\
$433$ & $11$ & $L_{2}(179^{2})$,\,$S_{4}(179)$,\,$U_{3}(199)$,\,$U_{4}(179)$\\
$439$ & $6$ & $L_{3}(439)$\\
$443$ & $9$ & $L_{2}(443^{2})$,\,$S_{4}(443)$\\
$449$ & $12$ & $L_{2}(67^{2})$,\,$S_{4}(67)$,\,$L_{4}(67)$\\
$457$ & $8$ & $L_{2}(109^{2})$,\,$S_{4}(109)$,\,$U_{4}(109)$\\
$463$ & $7$ & $L_{2}(463^{2})$,\,$S_{4}(463)$\\
$467$ & $17$ & $L_{2}(467^{2})$,\,$S_{4}(467)$,\,$U_{3}(467)$,\,$U_{4}(467)$\\
$487$ & $7$ & $U_{3}(233)$,\,$U_{4}(233)$\\
$491$ & $10$ & $U_{3}(491)$\\
$499$ & $11$ & $L_{3}(499)$,\,$L_{3}(139)$,\,$L_{3}(359)$,\,$L_{2}(499^{2})$,\,$S_{4}(499)$,\,$L_{4}(499)$\\
$509$ & $15$ & $L_{2}(509^{2})$,\,$S_{4}(509)$\\
$521$ & $19$ & $U_{5}(5)$,\,$L_{2}(5^{5})$,\,$U_{6}(5)$,\,$O^-_{10}(5)$,\,$L_{3}(521)$,\,$U_{3}(521)$,\,$U_{7}(5)$,\,$L_{5}(25)$,\\
 & &$S_{10}(5)$,\,$O_{11}(5)$,\,$O^+_{12}(5)$,\,$L_{2}(43^{4})$,\,$S_{4}(43^{2})$,\,$U_{8}(5)$,\,$L_{2}(521^{3})$,\,$G_2(521)$\\
$523$ & $24$ & $U_{3}(61)$,\,$U_{3}(463)$,\,$L_{2}(61^{3})$,\,$G_2(61)$,\,$U_{4}(463)$\\
$547$ & $18$ & $U_{3}(41)$,\,$U_{7}(3)$,\,$U_{4}(41)$,\,$U_{8}(3)$,\,$O^-_{14}(3)$,\,$U_{9}(3)$,\,$U_{10}(3)$\\
$557$ & $12$ & $L_{2}(439^{2})$,\,$S_{4}(439)$,\,$L_{2}(557^{2})$,\,$S_{4}(557)$,\,$L_{4}(439)$\\
$563$ & $8$ & $U_{3}(563)$\\
$571$ & $17$ & $L_{3}(109)$,\,$L_{4}(109)$,\,$L_{3}(461)$,\,$L_{2}(109^{3})$,\,$G_2(109)$,\,$L_{3}(571)$,\,$L_{3}(109^{2})$,\,$S_{6}(109)$,\\
 & &$O_{7}(109)$,\,$O^+_{8}(109)$\\
$577$ & $13$ & $L_{2}(577^{2})$,\,$S_{4}(577)$\\
$593$ & $11$ & $L_{2}(593^{2})$,\,$S_{4}(593)$,\,$U_{3}(593)$,\,$U_{4}(593)$\\
$599$ & $5$ & $L_{2}(599^{2})$,\,$S_{4}(599)$\\
$601$ & $27$ & $U_{3}(25)$,\,$^3D_4(5)$,\,$U_{4}(25)$,\,$L_{2}(5^{6})$,\,$S_{4}(5^{3})$,\,$G_2(25)$,\,$U_{3}(577)$,\,$L_{3}(5^{4})$,\\
 & &$S_{6}(25)$,\,$O_{7}(25)$,\,$F_4(5)$,\,$O^+_{8}(25)$,\,$L_{2}(601^{2})$,\,$S_{4}(601)$,\,$U_{4}(577)$,\,$O^-_{12}(5)$,\\
 & &$L_{6}(25)$,\,$O_{13}(5)$,\,$S_{12}(5)$,\,$O^-_{14}(5)$\\
$607$ & $16$ & $U_{3}(211)$,\,$U_{3}(397)$,\,$U_{4}(211)$,\,$L_{2}(211^{3})$,\,$G_2(211)$,\,$L_{3}(211^{2})$,\,$S_{6}(211)$,\,$O_{7}(211)$,\\
 & &$O^+_{8}(211)$\\
$613$ & $6$ & $L_{3}(547)$\\
$619$ & $15$ & $U_{3}(367)$,\,$U_{3}(619)$\\
$631$ & $25$ & $L_{3}(43)$,\,$L_{3}(587)$,\,$L_{4}(43)$,\,$L_{3}(631)$,\,$L_{2}(43^{3})$,\,$G_2(43)$,\,$L_{3}(43^{2})$,\,$S_{6}(43)$,\\
 & &$O_{7}(43)$,\,$O^+_{8}(43)$,\,$O^-_{8}(43)$,\,$L_{4}(43^{2})$,\,$O_{9}(43)$,\,$S_{8}(43)$\\
$641$ & $6$ & $L_{2}(487^{2})$,\,$S_{4}(487)$,\,$U_{3}(641)$\\
$653$ & $11$ & $L_{2}(149^{2})$,\,$S_{4}(149)$,\,$L_{3}(653)$,\,$L_{4}(149)$\\
$659$ & $5$ & $L_{2}(659^{2})$,\,$S_{4}(659)$\\
$661$ & $15$ & $U_{11}(3)$,\,$U_{12}(3)$\\
$673$ & $14$ & $U_{3}(2^{8})$,\,$^3D_4(16)$,\,$L_{2}(2^{24})$,\,$G_2(2^{8})$,\,$S_{4}(2^{12})$,\,$O^-_{8}(64)$,\,$L_{4}(2^{12})$,\,$S_{8}(64)$,\\
 & &$O^+_{10}(64)$\\
$677$ & $8$ & $U_{3}(677)$\\
$683$ & $19$ & $L_{2}(2^{11})$,\,$U_{11}(2)$,\,$U_{12}(2)$,\,$O^-_{22}(2)$,\,$L_{11}(4)$,\,$S_{22}(2)$,\,$O^+_{24}(2)$,\,$O^-_{24}(2)$,\\
 & &$L_{12}(4)$,\,$S_{24}(2)$\\
$701$ & $11$ & $L_{2}(701^{2})$,\,$S_{4}(701)$\\
$709$ & $16$ & $L_{3}(227)$,\,$L_{2}(613^{2})$,\,$S_{4}(613)$,\,$L_{2}(227^{3})$,\,$G_2(227)$\\
$719$ & $10$ & $U_{3}(719)$\\
$727$ & $8$ & $L_{3}(281)$\\
$733$ & $12$ & $L_{2}(353^{2})$,\,$S_{4}(353)$,\,$L_{3}(307)$,\,$L_{4}(307)$,\,$U_{4}(353)$\\
$739$ & $7$ & $U_{3}(419)$,\,$U_{3}(739)$\\
$743$ & $11$ & $L_{2}(743^{2})$,\,$S_{4}(743)$\\
$751$ & $9$ & $U_{3}(73)$,\,$U_{4}(73)$\\
$757$ & $17$ & $L_{3}(27)$,\,$L_{4}(27)$,\,$L_{2}(3^{9})$,\,$G_2(27)$,\,$L_{2}(757^{2})$,\,$S_{4}(757)$,\,$E_6(3)$,\,$L_{3}(3^{6})$,\\
 & &$S_{6}(27)$,\,$O_{7}(27)$,\,$O^+_{8}(27)$,\,$U_{6}(27)$\\
$769$ & $12$ & $U_{3}(19^{2})$,\,$^3D_4(19)$,\,$L_{2}(19^{6})$,\,$S_{4}(19^{3})$,\,$G_2(19^{2})$,\,$U_{3}(409)$,\,$U_{3}(3^{8})$\\
$773$ & $18$ & $L_{2}(317^{2})$,\,$S_{4}(317)$,\,$U_{3}(773)$\\
$787$ & $16$ & $L_{3}(787)$,\,$L_{3}(379)$,\,$L_{2}(787^{2})$,\,$S_{4}(787)$,\,$L_{4}(787)$\\
$809$ & $7$ & $L_{2}(491^{2})$,\,$S_{4}(491)$,\,$L_{3}(809)$,\,$U_{4}(491)$\\
$811$ & $16$ & $U_{3}(131)$,\,$L_{2}(811^{2})$,\,$S_{4}(811)$,\,$L_{3}(811)$,\,$L_{4}(811)$\\
$821$ & $7$ & $U_{3}(821)$,\,$L_{3}(821)$,\,$L_{2}(821^{3})$,\,$G_2(821)$\\
$823$ & $6$ & $L_{3}(823)$\\
$829$ & $18$ & $L_{3}(5^{3})$,\,$L_{4}(5^{3})$,\,$L_{2}(829^{2})$,\,$S_{4}(829)$,\,$U_{3}(829)$,\,$E_6(5)$,\,$U_{4}(829)$\\
$853$ & $7$ & $L_{2}(853^{2})$,\,$S_{4}(853)$\\
$857$ & $7$ & $L_{2}(857^{2})$,\,$S_{4}(857)$,\,$U_{3}(857)$,\,$U_{4}(857)$\\
$859$ & $8$ & $U_{3}(599)$,\,$U_{3}(859)$,\,$U_{4}(599)$\\
$863$ & $19$ & $U_{3}(863)$,\,$L_{2}(863^{2})$,\,$S_{4}(863)$,\,$U_{4}(863)$\\
$877$ & $10$ & $L_{2}(151^{2})$,\,$S_{4}(151)$,\,$U_{3}(283)$,\,$L_{2}(283^{3})$,\,$G_2(283)$\\
$881$ & $4$ & $U_{3}(881)$\\
$883$ & $7$ & $L_{3}(337)$,\,$L_{4}(337)$\\
$907$ & $6$ & $U_{3}(523)$\\
$911$ & $14$ & $L_{5}(19)$,\,$L_{6}(19)$,\,$U_{7}(7)$,\,$L_{2}(911^{2})$,\,$S_{4}(911)$\\
$919$ & $23$ & $U_{3}(53)$,\,$U_{4}(53)$,\,$L_{2}(53^{3})$,\,$G_2(53)$,\,$L_{3}(53^{2})$,\,$S_{6}(53)$,\,$O_{7}(53)$,\,$O^+_{8}(53)$,\\
 & &$L_{2}(919^{2})$,\,$S_{4}(919)$,\,$L_{3}(919)$,\,$L_{4}(919)$\\
$929$ & $10$ & $U_{3}(929)$\\
$937$ & $6$ & $U_{3}(13^{3})$\\
$941$ & $11$ & $L_{2}(97^{2})$,\,$S_{4}(97)$,\,$U_{4}(97)$,\,$U_{3}(941)$\\
$947$ & $8$ & $L_{3}(947)$\\
$953$ & $16$ & $U_{3}(953)$\\
$971$ & $8$ & $L_{3}(971)$\\
$991$ & $10$ & $L_{3}(113)$,\,$L_{3}(877)$,\,$L_{3}(991)$\\
$997$ & $14$ & $L_{3}(997)$\\
\hline
\end{longtable}
}

\bigskip
\bigskip
\bigskip

\footnotesize{
\begin{longtable}[h]{l|l|l}
\caption {\bf The simple groups $G$ with  $\pi(G)\se\{1,2,\ld,1000\}$}\label{alph}\\
\hline
$^{\vphantom{A^A}}$Series & Parameter & Values \\
\hline
$^{\vphantom{A^A}}$Sporadic& all & ---\\
$A_n$ & $n$ & $5,\ld,1008$\\
$L_{2}(p)$ & $p$ & $5,\ld,997$\\
$L_{2}(p^{2})$ & $p$ & $2,\ld,53$,\,$67$,\,$73$,\,$83$,\,$89$,\,$97$,\,$107$,\,$109$,\,$149$,\,$151$,\,$157$,\,$173$,\,$179$,\,$191$,\\
 & & $193$,\,$211$,\,$233,\ld,251$,\,$269$,\,$293,\ld,317$,\,$337$,\,$353$,\,$401$,\,$421$,\,$431$,\\
 & & $439$,\,$443$,\,$463$,\,$467$,\,$487$,\,$491$,\,$499$,\,$509$,\,$557$,\,$577$,\,$593$,\,$599$,\,$601$,\,$613$,\\
 & & $659$,\,$701$,\,$743$,\,$757$,\,$787$,\,$811$,\,$829$,\,$853$,\,$857$,\,$863$,\,$911$,\,$919$\\
$L_{2}(p^{3})$ & $p$ & $2,\ld,37$,\,$43$,\,$47$,\,$53$,\,$61$,\,$107$,\,$109$,\,$211$,\,$227$,\,$263$,\,$283$,\,$373$,\,$431$,\,$521$,\\
 & & $821$\\
$L_{2}(p^{k})$, $k\ge 4$ & $p^k$ & $2^{4}$,\,$2^{5}$,\,$2^{6}$,\,$2^{7}$,\,$2^{8}$,\,$2^{9}$,\,$2^{10}$,\,$2^{11}$,\,$2^{12}$,\,$2^{14}$,\,$2^{15}$,\,$2^{18}$,\,$2^{24}$,\,$3^{4}$,\,$3^{5}$,\,$3^{6}$,\,$3^{8}$,\\
 & & $3^{9}$,\,$5^{4}$,\,$5^{5}$,\,$5^{6}$,\,$7^{6}$,\,$19^{6}$,\,$29^{6}$,\,$43^{4}$\\
$L_{3}(p)$ & $p$ & $2,\ld,37$,\,$43$,\,$47$,\,$53$,\,$61$,\,$67$,\,$79$,\,$83$,\,$107$,\,$109$,\,$113$,\,$137$,\,$139$,\,$149$,\,$163$,\\
 & & $181$,\,$191$,\,$211$,\,$227$,\,$229$,\,$263$,\,$269$,\,$277$,\,$281$,\,$283$,\,$307$,\,$313$,\,$337$,\,$359$,\\
 & & $373$,\,$379$,\,$431$,\,$439$,\,$461$,\,$499$,\,$521$,\,$547$,\,$571$,\,$587$,\,$631$,\,$653$,\,$787$,\,$809$,\\
 & & $\ld,823$,\,$877$,\,$919$,\,$947$,\,$971$,\,$991$,\,$997$\\
$L_{3}(p^{2})$ & $p$ & $2,\ld,37$,\,$43$,\,$47$,\,$53$,\,$107$,\,$109$,\,$211$,\,$431$\\
$L_{3}(p^{k})$, $k\ge 3$ & $p^k$ & $2^{3}$,\,$2^{4}$,\,$2^{5}$,\,$2^{6}$,\,$2^{7}$,\,$2^{8}$,\,$2^{10}$,\,$2^{12}$,\,$3^{3}$,\,$3^{4}$,\,$3^{6}$,\,$5^{3}$,\,$5^{4}$\\
$L_{4}(p)$ & $p$ & $2,\ld,37$,\,$43$,\,$47$,\,$53$,\,$67$,\,$83$,\,$107$,\,$109$,\,$149$,\,$191$,\,$211$,\,$269$,\,$307$,\,$313$,\\
 & & $337$,\,$431$,\,$439$,\,$499$,\,$787$,\,$811$,\,$919$\\
$L_{4}(p^{k})$, $k\ge 2$ & $p^k$ & $2^{2}$,\,$2^{3}$,\,$2^{4}$,\,$2^{5}$,\,$2^{6}$,\,$2^{7}$,\,$2^{12}$,\,$3^{2}$,\,$3^{3}$,\,$3^{4}$,\,$5^{2}$,\,$5^{3}$,\,$43^{2}$\\
$L_n(q)$, $n=5,6$ & $q$ & $2$,\,$2^{2}$,\,$2^{3}$,\,$2^{4}$,\,$2^{6}$,\,$3$,\,$3^{2}$,\,$5$,\,$5^{2}$,\,$19$\\
$L_n(q)$, $n\ge 7$ & $n,q$ & $n=7$,\,$q=2,2^{2},2^{3},2^{4}$;\ $n=8$,\,$q=2,2^{2},2^{3}$;\ $n=9$,\,$q=2,2^{2}$;\\
 & & $n=10$,\,$q=2,2^{2}$;\ $n=11$,\,$q=2,2^{2}$;\ $n=12$,\,$q=2,2^{2}$\\
$U_{3}(p)$ & $p$ & $3,\ld,61$,\,$73$,\,$89,\ld,109$,\,$131$,\,$173$,\,$179$,\,$197$,\,$199$,\,$211$,\,$227$,\,$233$,\\
 & & $239$,\,$257$,\,$263$,\,$283$,\,$293$,\,$353$,\,$367$,\,$373$,\,$397,\ld,419$,\,$431$,\,$463$,\,$467$,\\
 & & $491$,\,$521$,\,$523$,\,$563$,\,$577$,\,$593$,\,$599$,\,$619$,\,$641$,\,$677$,\,$719$,\,$739$,\,$773$,\,$821$,\\
 & & $829$,\,$857$,\,$859$,\,$863$,\,$881$,\,$929$,\,$941$,\,$953$\\
$U_{3}(p^{k})$, $k\ge 2$ & $p^k$ & $2^{2}$,\,$2^{3}$,\,$2^{4}$,\,$2^{5}$,\,$2^{6}$,\,$2^{8}$,\,$3^{2}$,\,$3^{3}$,\,$3^{5}$,\,$3^{8}$,\,$5^{2}$,\,$7^{2}$,\,$13^{3}$,\,$19^{2}$,\,$29^{2}$\\
$U_{4}(p)$ & $p$ & $2,\ld,53$,\,$73$,\,$89$,\,$97$,\,$107$,\,$109$,\,$173$,\,$179$,\,$211$,\,$233$,\,$239$,\,$293$,\,$353$,\,$401$,\\
 & & $431$,\,$463$,\,$467$,\,$491$,\,$577$,\,$593$,\,$599$,\,$829$,\,$857$,\,$863$\\
$U_{4}(p^{k})$, $k\ge 2$ & $p^k$ & $2^{2}$,\,$2^{3}$,\,$2^{4}$,\,$2^{5}$,\,$2^{6}$,\,$3^{2}$,\,$3^{3}$,\,$5^{2}$\\
$U_n(q)$, $n=5,6$ & $q$ & $2$,\,$2^{2}$,\,$2^{3}$,\,$3$,\,$3^{3}$,\,$5$,\,$7$,\,$17$,\,$23$,\,$29$\\
$U_n(q)$, $n\ge 7$ & $n,q$ & $n=7$,\,$q=2,2^{2},3,5,7$;\ $n=8$,\,$q=2,2^{2},3,5$;\ $n=9$,\,$q=2,2^{2},3$;\\
 & & $n=10$,\,$q=2,2^{2},3$;\ $n=11$,\,$q=2,3$;\ $n=12$,\,$q=2,3$\\
$S_{4}(p)$ & $p$ & $3,\ld,53$,\,$67$,\,$73$,\,$83$,\,$89$,\,$97$,\,$107$,\,$109$,\,$149$,\,$151$,\,$157$,\,$173$,\,$179$,\,$191$,\\
 & & $193$,\,$211$,\,$233,\ld,251$,\,$269$,\,$293,\ld,317$,\,$337$,\,$353$,\,$401$,\,$421$,\,$431$,\\
 & & $439$,\,$443$,\,$463$,\,$467$,\,$487$,\,$491$,\,$499$,\,$509$,\,$557$,\,$577$,\,$593$,\,$599$,\,$601$,\,$613$,\\
 & & $659$,\,$701$,\,$743$,\,$757$,\,$787$,\,$811$,\,$829$,\,$853$,\,$857$,\,$863$,\,$911$,\,$919$\\
$S_{4}(p^{k})$, $k\ge 2$ & $p^k$ & $2^{2}$,\,$2^{3}$,\,$2^{4}$,\,$2^{5}$,\,$2^{6}$,\,$2^{7}$,\,$2^{9}$,\,$2^{12}$,\,$3^{2}$,\,$3^{3}$,\,$3^{4}$,\,$5^{2}$,\,$5^{3}$,\,$7^{3}$,\,$19^{3}$,\,$29^{3}$,\,$43^{2}$\\
$S_{6}(p)$ & $p$ & $2,\ld,37$,\,$43$,\,$47$,\,$53$,\,$107$,\,$109$,\,$211$,\,$431$\\
$S_{6}(p^{k})$, $k\ge 2$ & $p^k$ & $2^{2}$,\,$2^{3}$,\,$2^{4}$,\,$2^{5}$,\,$2^{6}$,\,$3^{2}$,\,$3^{3}$,\,$5^{2}$\\
$S_{8}(q)$ & $q$ & $2$,\,$2^{2}$,\,$2^{3}$,\,$2^{6}$,\,$3$,\,$3^{2}$,\,$5$,\,$43$\\
$S_n(q)$, $n\ge 10$ & $n,q$ & $n=10$,\,$q=2,2^{2},2^{3},3,5$;\ $n=12$,\,$q=2,2^{2},2^{3},3,5$;\\
 & & $n=14$,\,$q=2,2^{2}$;\ $n=16$,\,$q=2$;\ $n=18$,\,$q=2$;\\
 & & $n=20$,\,$q=2$;\ $n=22$,\,$q=2$;\ $n=24$,\,$q=2$\\
$O_{7}(p)$ & $p$ & $3,\ld,37$,\,$43$,\,$47$,\,$53$,\,$107$,\,$109$,\,$211$,\,$431$\\
$O_{7}(p^{k})$, $k\ge 2$ & $p^k$ & $3^{2}$,\,$3^{3}$,\,$5^{2}$\\
$O_n(q)$, $n\ge 9$ & $n,q$ & $n=9$,\,$q=3,3^{2},5,43$;\ $n=11$,\,$q=3,5$;\ $n=13$,\,$q=3,5$\\
$O^+_{8}(p)$ & $p$ & $2,\ld,37$,\,$43$,\,$47$,\,$53$,\,$107$,\,$109$,\,$211$,\,$431$\\
$O^+_{8}(p^{k})$, $k\ge 2$ & $p^k$ & $2^{2}$,\,$2^{3}$,\,$2^{4}$,\,$2^{5}$,\,$2^{6}$,\,$3^{2}$,\,$3^{3}$,\,$5^{2}$\\
$O^+_n(q)$, $n\ge 10$ & $n,q$ & $n=10$,\,$q=2,2^{2},2^{3},2^{6},3,3^{2},5$;\ $n=12$,\,$q=2,2^{2},2^{3},3,5$;\\
 & & $n=14$,\,$q=2,2^{2},2^{3}$;\ $n=16$,\,$q=2,2^{2}$;\ $n=18$,\,$q=2$;\\
 & & $n=20$,\,$q=2$;\ $n=22$,\,$q=2$;\ $n=24$,\,$q=2$\\
$O^-_{8}(q)$ & $q$ & $2$,\,$2^{2}$,\,$2^{3}$,\,$2^{6}$,\,$3$,\,$3^{2}$,\,$5$,\,$43$\\
$O^-_n(q)$, $n\ge 10$ & $n,q$ & $n=10$,\,$q=2,2^{2},2^{3},3,5$;\ $n=12$,\,$q=2,2^{2},2^{3},3,5$;\\
 & & $n=14$,\,$q=2,2^{2},3,5$;\ $n=16$,\,$q=2$;\ $n=18$,\,$q=2$;\\
 & & $n=20$,\,$q=2$;\ $n=22$,\,$q=2$;\ $n=24$,\,$q=2$\\
$G_2(p)$ & $p$ & $3,\ld,37$,\,$43$,\,$47$,\,$53$,\,$61$,\,$107$,\,$109$,\,$211$,\,$227$,\,$263$,\,$283$,\,$373$,\,$431$,\,$521$,\\
 & & $821$\\
$G_2(p^{k})$, $k\ge 2$ & $p^k$ & $2^{2}$,\,$2^{3}$,\,$2^{4}$,\,$2^{5}$,\,$2^{6}$,\,$2^{8}$,\,$3^{2}$,\,$3^{3}$,\,$5^{2}$,\,$7^{2}$,\,$19^{2}$,\,$29^{2}$\\
$F_4(q)$ & $q$ & $2$,\,$2^{2}$,\,$2^{3}$,\,$3$,\,$5$\\
$E_6(q)$ & $q$ & $2$,\,$2^{2}$,\,$3$,\,$5$\\
$E_7(q)$ & $q$ & $2$,\,$2^{2}$\\
$E_8(q)$ & $q$ & $2$\\
$^3D_4(q)$ & $q$ & $2$,\,$2^{2}$,\,$2^{3}$,\,$2^{4}$,\,$3$,\,$5$,\,$7$,\,$7^{2}$,\,$19$,\,$29$\\
$^2E_6(q)$ & $q$ & $2$,\,$2^{2}$,\,$3$\\
$Sz(2^{k})$ & $k$ & $3$,\,$5$,\,$7$,\,$9$\\
$^2G_2(3^{k})$ & $k$ & $3$,\,$5$\\
$^2F_4(2^{k})$, $k\ge 3$ & $k$ & $3$\\
$^2F_4(2)'$ & --- & ---\\
\hline
\end{longtable}
}

\end{document}